\documentclass[12pt,leqno]{article}
\usepackage{amsmath,amsthm,amssymb}
\usepackage[matrix,arrow]{xy}
\newcommand{\bP}{\mathbb P}
\newcommand{\A}{\mathbb A}

\newcommand{\cO}{{\mathcal O}}
\newcommand{\Spec}{\mathop{\rm Spec}\nolimits}
\newcommand{\Sym}{\mathop{\rm Sym}\nolimits}
\newcommand{\ord}{\mathop{\rm ord}\nolimits}
\newcommand{\Bl}{\mathop{{\rm B}\ell}\nolimits}
\newcommand{\smallcap}{\mathbin{\raise.4pt\hbox{$\scriptstyle\cap$}}}
\newcommand{\smallcirc}{\mathchoice
{\mathbin{\raise.6pt\hbox{$\scriptstyle\circ$}}}
{\mathbin{\raise.6pt\hbox{$\scriptstyle\circ$}}}
{\mathbin{\raise.35pt\hbox{$\scriptscriptstyle\circ$}}}
{\mathbin{\hbox{$\scriptscriptstyle\circ$}}}}
\newtheorem{thm}{Theorem}[section]
\newtheorem{pr}[thm]{Proposition}
\newtheorem{lm}[thm]{Lemma}
\newtheorem{cor}[thm]{Corollary}
\theoremstyle{definition}
\newtheorem{defn}[thm]{Definition}
\newtheorem*{keyobs}{Key Observation}
\newtheorem*{exas}{Examples}
\theoremstyle{remark}
\newtheorem{rem}[thm]{Remark}

\begin{document}
\title{Canonical rational equivalence of intersections of divisors}
\author{Andrew Kresch$^1$}
\date{3 December 1997}
\maketitle
\footnotetext[1]{Funded by a
Fannie and John Hertz Foundation Fellowship
for Graduate Study and an Alfred P. Sloan Foundation Dissertation
Fellowship.}

\setcounter{section}{-1}
\section{Introduction}

One way to define an operation in intersection theory is to define
a map on the group of algebraic cycles together with
a map on the group of rational equivalences
which commutes with the boundary operation.
Assuming the maps commute with smooth pullback, the extension of
the operation to the setting of algebraic stacks is automatic.
The goal of the first section of this paper
is to present the operation of intersecting
with a principal Cartier divisor in this light.

We then show how this operation lets us obtain a
rational equivalence which is fundamental to intersection theory.
A one-dimensional family of cycles on
an algebraic variety always admits a unique limiting cycle, but a family of
cycles over the punctured affine plane may yield different limiting
cycles if one approaches the origin from different directions.
An important step in the historical development of intersection theory was
realizing how to prove
that any two such limiting cycles are rationally equivalent.
The results of the first section yield, as a corollary,
a new, explicit formula for this rational equivalence.

Another important rational equivalence in intersection theory is
the one that is used to demonstrate commutativity of Gysin maps
associated to regularly embedded subschemes.
In section 2, we exhibit a two-dimensional family of cycles such
that the cycles we obtain from specializing in two different ways
are precisely the ones we need to show to be rationally equivalent to obtain
the commutativity result.
Our explicit rational equivalence
respects smooth pullback, and hence the generalization to stacks is
automatic.
This simplifies intersection theory on Deligne-Mumford stacks as in \cite{v},
where construction of such a rational equivalence
fills the most difficult section of that important paper.

Since our rational equivalence arises by considering families of cycles
on a larger total space,
we are able to deduce (section 3) that the rational
equivalence is invariant under a certain naturally arising group action.
The key observation is that we can manipulate the situation
so that the group action extends to the total space.
This equivariance result is used, but appears with mistaken proof,
in \cite{bf}, where an important new tool of modern intersection theory ---
the theory of virtual fundamental classes --- is developed.

The author would like to thank S. Bloch, W. Fulton, T. Graber,
and R. Pandharipande for helpful advice and
the organizers and staff of the Mittag-Leffler Institute for
hospitality during the 1996--97 program in algebraic geometry.

\section{Intersection with divisors}

In this section we work exclusively on schemes of finite type
over a fixed base field.
The term variety denotes integral scheme, and by a subvariety we mean an
integral closed subscheme.
We denote by $Z_*X$, $W_*X$, and $A_*X$, respectively, the group of
algebraic cycles, group of rational equivalences, and Chow group
of a scheme $X$.
The boundary map $W_*X\to Z_*X$ is denoted $\partial$.
We refer to \cite{f} for basic definitions and properties from
intersection theory.
Given a Cartier divisor $D$ we denote by $[D]$ the associated
Weil divisor
(it is important to note that the notion of Weil divisor makes
sense on arbitrary varieties, \cite{f} \S 1.2).
If $X$ is a variety then we
denote by $X^1$ the set of subvarieties of
codimension 1.

\begin{defn}
Let $X$ be a variety and let $D$ be a Cartier divisor.
Let $\pi\colon \widehat X\to X$
be the normalization map.
The {\em support} of $D$, denoted $|D|$, is defined to be
$\pi(\bigcup_{\substack{W\in\widehat X^1\\ \ord_W \pi^*D\ne 0}} W).$
\end{defn}

\begin{rem}
This agrees with the na\"\i{}ve notion of support (the union of
all subvarieties appearing with nonzero coefficient in
$[D]$) when $X$ is normal or when $D$ is effective.
\end{rem}

\begin{rem}
There is yet another notion of support which appears in \cite{f}.
There, the support of a divisor is a piece of data
that must be specified along with the divisor.
Given a Cartier divisor $D$ on a variety $X$, let $Z$ be any
closed subscheme such that away from $Z$ the canonical section
of $\cO(D)$ is well-defined and nonvanishing.
Then, \cite{f} defines an intersection operation
$A_k(X)\to A_{k-1}(Z)$.
Unfortunately, the support $|D|$ which we have defined
is not generally a support in this sense.
Hence in the definition below we require that our divisors be
specified by defining functions which are regular away from
their supports.
\end{rem}

We shall denote by $|D|^0$ the set of irreducible components of $|D|$.

\begin{defn}
\label{pdivisor}
Let $X$ be a variety.
A {\em $P$-divisor} on $X$ is a tuple $(U,U',x)$ such that
\begin{itemize}
\vspace{-12pt}
\item[(i)] $U$ and $U'$ are nonempty open subschemes of $X$ such
that $U\cup U'=X$;
\item[(ii)] $x\in k(U)^*$;
\item[(iii)] $x|_{U\cap U'}\in \cO^*(U\cap U')$; and
\item[(iv)] the data $(x\in k(U)^*, 1\in k(U')^*)$ specifies
a Cartier divisor $D$ such that $|D|=X\setminus U'$.
\end{itemize}
\end{defn}

By abuse of terminology, we call $D$ a $P$-divisor
if $D$ is the Cartier divisor associated to a
$P$-divisor as in (iv).
Given a $P$-divisor as above, we call
$x$ the {\em local defining function}.
A $P$-divisor may be pulled back via a morphism of varieties
provided that the image of the morphism is not
contained in the support of the underlying Cartier divisor.

\begin{exas}
\begin{itemize}
\item[] \hspace{-30pt} (i)\hspace{5pt}Let $X$ be a normal variety.
Let $x\in k(X)^*$ specify a principal Cartier divisor $D$.
Then $(X, X\setminus |D|, x)$ is a $P$-divisor.
\item[(ii)] Let $X$ be a variety.
Every effective principal Cartier divisor is a $P$-divisor.
\item[(iii)] Let $X$ be a variety, and let $\pi\colon X\to \bP^1$ be
a dominant morphism.
Then the fiber of $\pi$ over $\{0\}$ is a $P$-divisor.
\end{itemize}
\end{exas}

The operation of intersecting with a Cartier divisor is
generally defined only on the level of rational equivalence classes of cycles.
When $V\subset |D|$, we have
$D\cdot[V]=c_1(\cO(D)|_V)\smallcap[V]$, and
there is generally no way to pick canonically
a cycle representing this first Chern class.
The exception is when $\cO(D)|_{|D|}$ is trivial, or in our terminology,
$D$ is a $P$-divisor.
Then, we may define a cycle-level intersection operation
(see \cite{f}, Remark 2.3).

\begin{defn}
Let $X$ be a variety, and let $D$ be a $P$-divisor on $X$.
The cycle-level intersection operation
$$D\cdot {-}\colon Z_k(X)\to Z_{k-1}(|D|)$$
is given by
$$D\cdot [V] = \begin{cases}
{}[D|_V] & \text{if $V\not\subset |D|$}; \\
0& \text{if $V\subset |D|$}.
\end{cases}$$
\end{defn}

The claim that this map passes to rational equivalence and hence gives an
intersection operation
$D\cdot{}\colon A_k(X)\to A_{k-1}(|D|)$
is proved in \cite{f}, but not in a way that makes it easy to see
how $D\cdot\alpha$ is to be rationally equivalent to zero if
$\alpha$ is a cycle that is rationally equivalent to zero.
Following the program set out in the introduction, we would like
to demonstrate this fact by giving an explicit map on rational equivalences
which commutes with the boundary operation.

\begin{defn}
Let $X$ be a variety, and let $D$ be a $P$-divisor on $X$ with
local defining function $x$.
Say $V$ is a subvariety of $X$ with normalization
$\pi\colon \widehat V\to V$,
and suppose $y\in k(V)^*$.
We define the intersection operation on the level of rational equivalences
$$D\cdot{-}\colon W_k(X)\to W_{k-1}(|D|)$$
by
\begin{equation}
\label{maponrat}
D\cdot y = \begin{cases}
{}\pi_*\bigl(\bigoplus_{W\in |\pi^*D|^0} (y^{\ord_W x} / x^{\ord_W y}) |_W
\bigr) &
\text{if $V\not\subset |D|$}; \\
0& \text{if $V\subset |D|$}.
\end{cases}
\end{equation}
Here, $\pi_*\colon W_*\widehat V\to W_*V$ is pushforward of
rational equivalence.
\end{defn}

\begin{rem}
This definition explains why we a require the definition of a $P$-divisor
to include more
data than just that of the underlying Cartier divisor.
The map (\ref{maponrat}) actually depends on the choice of defining function.
\end{rem}

\begin{pr}
Let $X$ be a variety and let $D$ be a $P$-divisor on $X$.
Then the diagram
$$
\xymatrix{
W_k(X) \ar[r]^(.43){D\cdot{}}\ar[d]_\partial & W_{k-1}(|D|) \ar[d]_\partial \\
Z_k(X) \ar[r]^(.43){D\cdot{}} & Z_{k-1}(|D|)
}
$$
commutes.
\end{pr}

This follows easily from
\begin{pr}
Let $X$ be a normal variety and let $x$ and $y$ be rational functions
with associated principal Cartier divisors $D$ and $E$.
For $V\in X^1$ set $a_V=\ord_V x$ and $b_V=\ord_V y$.
Then
\begin{align}
\sum_{V\in X^1} \partial(y^{a_V}/x^{b_V}|_V)&=0;  \label{eqone} \\
\partial(D\cdot y) &= D\cdot(\partial\, y); \label{eqtwo} \\
D\cdot[E] - E\cdot[D] &= \sum_{V\in |D|^0\cap |E|^0}
\partial(y^{a_V}/ x^{b_V}|_V). \label{eqthree}
\end{align}
\end{pr}

\begin{proof}
If we split the sum in (\ref{eqone}) into a sum over $V\in |D|^0$ and
a sum over $V\not\in |D|^0$ we obtain (\ref{eqtwo}).
Similarly if we split away the terms with $V\in |D|^0\cap |E|^0$
we obtain (\ref{eqthree}) from (\ref{eqone}).
So, for a fixed variety $X$ and fixed divisors $D$ and $E$, the
three assertions are equivalent.
Now, we get (\ref{eqone}) as a consequence of the tame symbol in $K$-theory,
cf.\ \cite{q} \S7,
or by the following elementary geometric argument.
We quickly reduce to the case where
$D$ and $E$ are effective.
Then, when $D$ and $E$ meet properly,
(\ref{eqthree}) follows from \cite{f}, Theorem 2.4, case 1.
An induction on {\em excess of intersection}
$$\varepsilon(D,E)=\max_{V\in X^1} a_V\cdot b_V$$
completes the proof: if we denote the normalized
blow-up along the ideal $(x,y)$ by $\sigma\colon X'\to X$
and denote the exceptional divisor by $Z$ then we may write
$\sigma^*D = Z + D'$ and $\sigma^*E = Z + E'$,
and now $|D'|\cap |E'|=\emptyset$
and
$\max(\varepsilon(D',Z), \varepsilon(E',Z))<
\varepsilon(D,E)$
(assuming $D$ and $E$ do not meet properly), cf.\ \cite{f}, Lemma 2.4.
The result pushes forward.
\end{proof}

\begin{cor}
\label{canrat}
Let $D$ and $E$ be $P$-divisors on a variety $X$,
with respective local defining functions $x$ and $y$.
Let $\pi\colon\widehat X\to X$ be the normalization map.
Then
$$D\cdot[E] - E\cdot[D] = \partial\, \omega$$
where $\omega\in W_*(|D|\cap |E|)$ is given by
$$\omega = \sum_{V\in |\pi^*D|^0\cap |\pi^*E|^0}
\pi_*(y^{\ord_Vx}/ x^{\ord_Vy}|_V).$$
\end{cor}

\section{Application to intersection theory on stacks}

All stacks (and schemes) in this section are algebraic stacks of Artin type,
\cite{a}, \cite{l}, which are locally of finite type over the base field.
The notion of $P$-divsor on a stack makes sense (it is as in
Definition \ref{pdivisor} with ``open subscheme'' replaced by
``open substack,''
where by ``Cartier divisor'' in part (iv) of the definition we
mean a global section of the sheaf ${\mathcal K}^*/\cO^*$ for the
Zariski topology, and where
normalization, order along a substack of codimension 1,
and support of a Cartier divisor are well defined on stacks because
they all respect smooth pullback and hence can be defined locally).
Since an Artin stack possesses
a smooth cover by a scheme, the operation of intersecting
with a $P$-divisor on a stack comes for free
once we know that this operation on schemes commutes with
smooth pullback.
Also for free we get Corollary \ref{canrat} in the setting of stacks:
the formation of $\omega$ from $X$, $D$, and $E$ commutes with smooth
pullback.

\begin{pr}
Let $X$ be a variety, let $Y$ be a scheme,
and let $f\colon Y\to X$ be a smooth morphism.
Let $D$ be a $P$-divisor on $X$.
Then $f^*\smallcirc D\cdot{} = (f^*D)\cdot{}\smallcirc f^*$,
both as maps on cycles and as maps on rational equivalences.
\end{pr}

We now turn to an application of Corollary \ref{canrat} to
intersection theory on Deligne-Mumford stacks
(where a reasonable intersection theory exists, cf.\ \cite{g}, \cite{v}).
Central to intersection theory on schemes is the Gysin map corresponding
to a regularly embedded subscheme, since the diagonal of
a smooth scheme is a regular embedding and this way we obtain
an intersection product on smooth schemes.
The diagonal morphism for a smooth Deligne-Mumford stack is not generally
an embedding, but it is representable and unramified.

\begin{lm}
\label{unram}
Let $f\colon F\to G$ be a representable morphism of Artin stacks.
Then $f$ is unramified if and only if there exists
a commutative diagram
$$
\xymatrix{
U\ar[r]^g\ar[d] & V \ar[d] \\
F\ar[r]^f & G
}
$$
such that the vertical maps are smooth surjective,
$g$ is a closed immersion of schemes,
and the induced morphism $U\to F\times_GV$ is \'etale.
\end{lm}

\begin{proof}
This is \cite{v}, Lemma 1.19.
Because this is such a basic fact about properties of
morphisms in algebraic geometry, we present an elementary proof
in the Appendix.
\end{proof}

To describe a representable morphism,
we use the terminology {\em local immersion} as a synonym for
{\em unramified} and call $f$ above a
{\em regular local immersion} if moreover $g$ is a regular
embedding of schemes.
Since formation of normal cone is of a local nature, an obvious
patching construction produces the normal cone $C_XY$ to a local
immersion $X\to Y$; the cone is a bundle in case
$X\to Y$ is a regular local immersion.

To get Fulton-MacPherson-style intersection theory on Deligne-Mumford stacks
we clearly need to have Gysin maps for regular local immersions.
In \cite{v}, the author supplies this needed Gysin map by giving
a (long, difficult) proof of
the stack analogue of \cite{f}, Theorem 6.4, namely

\begin{pr}
\label{bigrat}
Let $X\to Y$ and $Y'\to Y$ be local immersions of Artin stacks.
Then
$[C_{X\times_YC_{Y'}Y}C_{Y'}Y]=[C_{C_XY\times_YY'}C_XY]$
in $A_*(C_XY\times_YC_{Y'}Y)$.
\end{pr}

\begin{rem}
Though our focus is on applications to intersection theory on
Deligne-Mumford stacks, we continue to make use of constructions which
behave well locally with respect to smooth pullback, and hence
our results are valid in the more general setting of Artin stacks.
\end{rem}

\begin{rem}
Given a stack $X$ which is only locally of finite type over a base field,
we must take $Z_*X$ to be the group of {\em locally finite} formal
linear combinations of integral closed substacks.
More intrinsically, $Z_*X$ is the group of global sections of the
sheaf for the smooth topology $\mathcal Z_*$
which associates to a stack of finite type
the free abelian group on integral closed substacks.
Similarly, $W_*X$ is the group of global sections of sheaf $\mathcal W_*$.
As always, $A_*X$ is defined to be $Z_*X/\partial W_*X$.
\end{rem}

The methods of the last section allow us to supply a new, simpler proof of
this proposition.

\begin{proof}
Recall that given a closed immersion $X\to Y$ there are associated
spaces
\begin{align*}
M_XY &= \Bl_{X\times\{0\}}Y\times\bP^1, \\
M^\circ_XY &= M_XY \setminus \Bl_{X\times\{0\}}Y\times\{0\},
\end{align*}
cf.\ \cite{f} \S 5.1.
Given a locally closed immersion, say with
$U$ is an open subscheme of $Y$
and $X$ a closed subscheme of $U$, then
$M^\circ_XY:=M^\circ_XU\amalg_{U\times\A^1} Y\times\A^1$ makes sense
and is independent of the choice of $U$.

This lets us define $M^\circ_FG$ when $F\to G$ is a local
immersion of stacks, as follows.
Assume we have a diagram as in the statement of Lemma \ref{unram},
and set $R=U\times_FU$ and $S=V\times_GV$.
There are projections $q_1, q_2\colon S\to G$.
Define $s_i\colon M^\circ_RS\to M^\circ_UV$ ($i=1,2$) to be the composite
$M^\circ_RS\to M^\circ_{U\times_GV}S\to M^\circ_UV$,
where the first map is induced by the open immersion $R\to U\times_GV$
and the second, by pullback via $q_i$.
Then $[M^\circ_RS\rightrightarrows M^\circ_UV]$ is the smooth
groupoid presentation of a stack which we denote $M^\circ_FG$.
We have, by descent, a morphism $M^\circ_FG\to\bP^1$, which is flat and
has as general fiber a copy of $G$ and as special fiber the
normal cone $C_FG$.

In the situation at hand, this construction gives
$$(s\times t)\colon M^\circ_XY\times_Y M^\circ_{Y'}Y\to \bP^1\times\bP^1,$$
and hence a pair of $P$-divisors,
$D$ (corresponding to $s$)
and $E$ (corresponding to $t$).
We note that
$(s\times t)^{-1}(\{0\}\times\{0\})=C_XY\times_YC_{Y'}Y$.
Since the restriction of $s\times t$ to
$\bP^1\times\bP^1\setminus \{0\}\times\{0\}$ is flat, we have
\begin{align*}
{}[D] &= [C_XY\times_Y M^\circ_{Y'}Y] \mod Z_*(C_XY\times_YC_{Y'}Y), \\
{}[E] &= [M^\circ_XY\times_Y C_{Y'}Y] \mod Z_*(C_XY\times_YC_{Y'}Y).
\end{align*}

We examine the fiber of $s\times t$ over $\bP^1\times\{0\}$ more closely.
The fiber square
$$\xymatrix{
i^*C_{Y'}Y \ar[r] \ar[d] & C_{Y'}Y \ar[d] \\
X \ar[r]^i & Y
}$$
gives rise to a closed immersion $f$ making
$$\xymatrix@C=2pt{
M^\circ_{i^*C_{Y'}Y}C_{Y'}Y \ar[rr]^(.48)f\ar[dr]_h &&
M^\circ_XY \times_Y C_{Y'}Y \ar[dl]^g \\
& {\bP^1}
}$$
commute (where $g$ is first projection followed by $s$).
Since $f$ is an isomorphism away from the fiber over 0,
we see in fact that
$$[E] = [M^\circ_{i^*C_{Y'}Y}C_{Y'}Y] \mod Z_*(C_XY\times_YC_{Y'}Y),$$
and since $h$ is flat we find
$$D\cdot [E] = [C_{i^*C_{Y'}Y}C_{Y'}Y].$$
Similarly, if $j$ denotes the map $Y'\to Y$ then
$$E\cdot [D] = [C_{j^*C_XY}C_XY]$$
and so the rational equivalence $\omega\in W_*(C_XY\times_YC_{Y'}Y)$ of
Corollary \ref{canrat} satisfies
$$\partial\,\omega = [C_{X\times_YC_{Y'}Y}C_{Y'}Y] -
[C_{C_XY\times_YY'}C_XY]. \qed $$
\renewcommand{\qed}{}\end{proof}

\begin{rem}
The map $M^\circ_FG\to G$ associated to a local immersion of
stacks is not generally separated, though this should cause
the reader no concern, since intersection theory is valid
even on non-separated schemes and stacks.
In fact, even those operations of \cite{v} which require a
so-called finite parametrization may be carried out
on arbitrary Deligne-Mumford stacks which are of finite type
over a field (no such operations show up in this paper).
This is so thanks to the proof, \cite{l} (10.1), that
every Deligne-Mumford stack of finite type over a field
possesses a finite parametrization, i.e.,
admits a finite surjective map from a scheme.
\end{rem}

\begin{rem}
The reader who wishes greater generality may see easily that
all results in this section are valid in the setting of
Artin stacks which are locally of finite type over an
excellent Dedekind domain.
\end{rem}

\section{Equivariance for tangent bundle action}

We continue to work with stacks which are locally of finite type
over some base field.
A special case of Proposition \ref{bigrat} is when
$i\colon X\to Y$ is a local immersion of smooth Deligne-Mumford stacks.
Suppose $j\colon Y'\to Y$ is a local immersion,
with $Y'$ an arbitrary Deligne-Mumford stack.
Recall that the local immersion $j$ gives rise
to a natural group action of $j^*T_Y$ on $C_{Y'}Y$.
In short, the action is given locally (say $Y$ is an affine scheme and
$Y'$ is the closed subscheme given by the ideal $I$) by
considering the action of $T_Y|_{Y'}$ on $\Spec \Sym (I/I^2)$
induced by the map $I/I^2\to \Omega^1_Y$ and proving
(\cite{bf}, Lemma 3.2) that the
normal cone $\Spec \bigoplus I^k/I^{k+1}$ is invariant under
the group action.

If we let $N_XY$ be the normal bundle to $X$ in $Y$ and denote simply
by $N$ its pullback to $X':=X\times_YY'$,
then $C_XY\times_YC_{Y'}Y$ is identified with $N\times_{X'} i^*C_{Y'}Y$.
Viewing $T_{Y'}$ as a subbundle of $j^*T_Y$,
we have the natural action of $T_{Y'}|_{X'}$ on
$i^*C_{Y'}Y$.
This plus the trivial action on $N$ gives an action
of $T_{Y'}|_{X'}$ on $N\times_{X'}i^*C_{Y'}Y$.

\begin{thm}
The rational equivalence between
$[C_{i^*C_{Y'}Y}C_{Y'}Y]$ and $[N\times_{X'}C_{X'}Y']$
produced in the proof of Proposition \ref{bigrat}
is invariant under the action of
$T_{Y'}|_{X'}$ on $N\times_{X'}C_{X'}Y'$
described above.
\end{thm}

As a consequence, the rational equivalence
descends to a rational equivalence on
the stack quotient $[N\times_{X'}i^*C_{Y'}Y\,/\,T_{Y'}|_{X'}]$.
This fact is exploited in \cite{bf}
(Lemma 5.9, where the authors invoke the incorrect stronger
claim appearing in Proposition 3.5 that the rational equivalence
is equivariant for the bigger group $T_Y|_{X'}$).

\begin{proof}
The question is local, so we may assume $Y$ is an irreducible
scheme, smooth and of finite type over the base field,
$X$ is an smooth irreducible closed subscheme of $Y$,
and $Y'$ is a closed subscheme of $Y$.
If $X\subset Y'$ then the group action is trivial and there is
nothing to prove, so we assume the contrary.

\begin{lm}
Let $Y$ be a smooth irreducible scheme of finite type over
a field $k$, of dimension $n$, let $X$ be a
smooth irreducible closed subscheme of $Y$ of codimension $d$, and
let $Y'$ be a closed subscheme of $Y$ such that $X\not\subset Y'$.
Let $x$ be a closed point of $Y'\cap X$.
Then, after suitable base change by a finite separable extension of
the base field, and after shrinking $Y$ to a neighborhood of $x$ in $Y$,
there exists an \'etale map $f\colon Y\to \A^n$
such that $X$ maps into a linear subspace of $\A^n$ of codimension $d$
and such that $Y'\to f(Y')$ is \'etale.
\end{lm}

\begin{proof}
We may assume $x$ is a $k$-valued point,
and moreover that $Y$ sits in $\A^l$ with $X=\A^{l-d}\cap Y$
(for suitable $l$).
We may take $x$ to be the origin of $\A^l$.
We consider as candidates for $f$ all linear functions mapping
the flag $\A^{l-d}\subset \A^l$ into the flag $\A^{n-d}\subset \A^n$.

Those $f$ with $f_*\colon T_{x,Y}\to T_{f(x),\A^n}$ surjective
form an open subscheme $U$ of $\A^{nl-dl+d^2}$.
Define locally closed subschemes $V_1$ and $V_2$ of $Y\times U$ by
$$V_1=\{(y,f)\in (Y'\cap X\setminus\{x\})\times U\,|\,
f(y)=0\}$$
and
$$V_2=\{(y,f)\in (Y'\setminus X)\times U\,|\,f(y)=0\},$$
and let $pr_2\colon Y\times U\to U$ be projection.
A dimension count using the fact that $X\not\subset Y'$
gives $\dim(V_1) < \dim(U)$ and $\dim(V_2) < \dim(U)$,
and hence
$U\setminus \bigl( \overline{pr_2(V_1)}\cup \overline{pr_2(V_2)} \bigr)$
is nonempty.
\end{proof}

Since the rational equivalence of the proof of
Proposition \ref{bigrat} commutes with \'etale base change,
we are reduced by the Lemma to the case where
$Y=\A^n$ and $X=\A^m$ (as a linear subspace of $\A^n$).
Now we need the

\begin{keyobs}
Assume $Y=\A^n$ and $Y'$ is a closed subscheme of $Y$.
Identify $T_Y$, as a group scheme over $Y$,
with the additive group $\A^n$.
Then there is a group action of $\A^n$ on $\widetilde M^\circ_{Y'}Y$
(which we define to be the fiber of $M^\circ_{Y'}Y\to\bP^1$ over $\A^1$)
which restricts to the natural action of $T_Y$
on $C_{Y'}Y$.
\end{keyobs}

Indeed, we let $\A^n$ act on $Y\times\A^1$ by
$$(a_1,\ldots,a_n)\cdot(x_1,\ldots,x_n,t)=(x_1+ta_1,\ldots,x_n+ta_n).$$
By the universal property of blowing up, this extends uniquely to an
action of $\A^n$ on $\widetilde M^\circ_{Y'}Y$.
If $Y'$ is given by the ideal $(f_1,\ldots,f_k)$, and if
we view $\widetilde M^\circ_{Y'}Y$ as the closure of
the graph of
$(f_1/t,\ldots,f_k/t)\colon Y\times(\A^1\setminus\{0\})\to
\A^k=\Spec k[z_1,\ldots,z_k]$, then
the action is given coordinatewise by
$${\mathbf a}=(a_1,\ldots,a_n)\colon z_i\mapsto z_i+
(f_i({\mathbf x} + t{\mathbf a}) - f_i({\mathbf x}))/t,$$
so at $t=0$ we recover $z_i\mapsto z_i+D_{\mathbf a}f_i({\mathbf x})$.
This is the natural action of $T_Y$ on $C_{Y'}Y$.

Concluding the proof of equivariance, we observe that
$\widetilde M^\circ_{\A^m}\A^n$ fits into the fiber diagram
$$
\xymatrix{
{\widetilde M^\circ_{\A^m}\A^n\times_{\A^n}\widetilde M^\circ_{Y'}\A^n}
\ar[r] \ar[d] & {\widetilde M^\circ_{Y'}\A^n} \ar[d] \\
{\widetilde M^\circ_{\A^m}\A^n} \ar[r] \ar[d] &
{\A^n} \ar[d] \\
{\widetilde M^\circ_{\{0\}}\A^{n-m}} \ar[r] &
{\A^{n-m}}
}$$
and now the action from the Key Observation of $\A^m\subset \A^n$ on
$\widetilde M^\circ_{Y'}Y$, plus the trivial action of $\A^m$ on
$M^\circ_{\{0\}}\A^{n-m}$, combine to give a group action
of $\A^m$ on $\widetilde M^\circ_XY\times_Y\widetilde M^\circ_{Y'}Y$.
The function
$\widetilde M^\circ_XY\times_Y\widetilde M^\circ_{Y'}Y\to \A^1\times\A^1$
which is used in Corollary \ref{canrat} is invariant for
this $\A^m$-action.
Since the rational equivalence of the proof of
Proposition \ref{bigrat} is compatible with smooth pullback,
we get the desired equivariance result.
\end{proof}

\section{Appendix: unramified morphisms}

We give an elementary algebraic proof of the following fact.

\begin{lm}
Let $S\to T$ be an unramified morphism of affine schemes which are of
finite type over a base field $k$.
Then there exists a commutative diagram of affine schemes
$$\xymatrix{
U \ar[r]^g \ar[d] & V \ar[d] \\
S \ar[r]^f & T
}$$
such that the vertical maps are \'etale surjective
and such that $g$ is a closed immersion.
\end{lm}

This fact plus the local nature of the property of being unramified
gives us Lemma \ref{unram}.

\begin{proof}
Say $S=\Spec A$, $T=\Spec B$, and $f$ is given algebraically
by $f^*\colon B\to A$.
Recall that for $f$ to be unramified means that
for every maximal ideal ${\mathfrak p}$ of $A$ with
${\mathfrak q}=f({\mathfrak p})$, we have
$f^*({\mathfrak q})\cdot A_{\mathfrak p}={\mathfrak p}A_{\mathfrak p}$,
and the induced field extension
$B/{\mathfrak q}\to A/{\mathfrak p}$ is separable.

{\em Case 1:} The induced field extension
$B/{\mathfrak q}\to A/{\mathfrak p}$ is an isomorphism.
Then, if $x_1,\ldots,x_n$ are generators of $A$ as a $k$-algebra,
we may write
$$x_i=f^*(t_i)+w_i$$
with $t_i\in B$ and $w_i\in {\mathfrak p}$, for each $i$.
Since $f$ is unramified, we have
$$w_i=\sum_{j=1}^{m_i} \frac{f^*(y_{ij})p_{ij}}{q_i}$$
for some $y_{ij}\in{\mathfrak q}$, $p_{ij}\in A$, and
$q_i\in A\setminus{\mathfrak p}$.

Choose representative polynomials $P_{ij}$ and $Q_i$ in
$k[X_1,\ldots,X_n]$ such that
$P_{ij}(x_1,\ldots,x_n)=p_{ij}$ and
$Q_i(x_1,\ldots,x_n)=q_i$.
Let
\begin{eqnarray*}
\lefteqn{
V=\Spec B[X_1,\ldots, X_n]\bigm/\big(\,X_1Q_1-t_1Q_1-\sum_{j=1}^{m_1}
y_{1j}P_{1j}\,,\ \ldots,}\hspace{130pt} \\
& & X_nQ_n-t_nQ_n-\sum_{j=1}^{m_n}y_{nj}P_{nj}\,\big),
\end{eqnarray*}
and define
$g\colon S\to V$ by $B\stackrel{f^*}\rightarrow A$ and $X_i\mapsto x_i$,
and let $\varphi\colon V\to T$ be given by inclusion of $B$.
Then $g$ is a closed immersion, and by the Jacobian criterion
$\varphi$ is \'etale in some neighborhood of $g({\mathfrak p})$.

{\em Case 2:} The field extension
$B/{\mathfrak q}\to A/{\mathfrak p}$ is separable.
Let $k'$ be the maximal subfield of $A/{\mathfrak p}$ which is
separable over $k$,
and make the \'etale base change
$\Spec k'\to \Spec k$ to get
$f'\colon S'\to T'$.
Now $S'$ has an $A/{\mathfrak p}$-valued point which maps to
${\mathfrak p}\in S$, and since $k'$ together with
$B/{\mathfrak q}$ generates all of $A/{\mathfrak p}$
we are now in the situation of Case 1.
\end{proof}

\noindent
Department of Mathematics\\
University of Chicago\\
5734 S. University Avenue\\
Chicago, IL 60637\\
kresch@math.uchicago.edu
\end{document}